\documentclass{amsart}
\usepackage{amsfonts,amssymb}
\usepackage{epsf}

\newcommand{\R}{\mathbb{R}}
\newcommand{\N}{\mathbb{N}}
\newcommand{\Z}{\mathbb{Z}}
\newcommand{\Id}{\mathrm{Id}}

\newtheorem{defn}{Definition}
\newtheorem{thm}{Theorem}
\newtheorem{lem}{Lemma}

\DeclareMathOperator{\supp}{supp}

\begin{document}

\title[Multiple Recurrence and Tiling]{An Application of Topological
  Multiple Recurrence to Tiling}

\author{R. de la Llave} \address{Mathematics Department\\University of
  Texas at Austin\\1 University Station, C1200\\TX 78712-0257\\U.S.A.}
\author{A. Windsor} \address{Mathematical Sciences\\373 Dunn
  Hall\\University of Memphis\\TN 38152-3240\\U.S.A.}

\begin{abstract} 
  We show that given any tiling of Euclidean space, any geometric
  patterns of points, we can find a patch of tiles (of arbitrarily
  large size) so that copies of this patch appear in the tiling nearly
  centered on a scaled and translated version of the pattern.  The rather
  simple proof uses 
  Furstenberg's topological multiple recurrence theorem.


\end{abstract} 

\maketitle

\section{Introduction}
Tilings have a very long story as an applied sub-discipline at the
service of architecture and decoration~\cite{Grunbaum}.

The modern mathematical theory of tiling (we collect several of the
standard definitions of tilings in Section~\ref{sec:definitions}) can
be said to begin in the 60's with the work \cite{Wang}. The surprising
main result of \cite{Wang} is that there is no algorithmic way to
prove that certain collections of tiles in $\R^d$ could be used to
tile the whole of $\R^d$, see also \cite{Berger}. Central to the proof
is the construction of sets of tiles that tile the plane only
\emph{non-periodically}, such sets of tiles are called
\emph{aperiodic}.  By now many aperiodic sets of tiles are known of
varying complexity~\cite{Grunbaum}. Perhaps the most famous aperiodic
set of tiles are the Penrose tiles~\cite{Penrose}.


The goal of this paper is sort of complementary with the previous
results.  Rather than showing that there are aperiodic sets of tiles,
or sets of tiles that exhibit unordered behavior~\cite{Quaquaversal},
we show that for any tiling, there has to be some amount of order, in
particular approximate periodicity. Roughly speaking, our results show
(see Theorems~\ref{MainThm1}, \ref{MainThm2}, \ref{MainThm3} for
precise formulations) that in \emph{any} tiling, given \emph{any}
pattern and \emph{any} size $R$, we can find a patch of size at least
$R$ so that copies of this patch appear in the tiling nearly centered
on a scaled and translated copy of the specified pattern.  We can also
require that the scalings of the pattern appearing in the conclusions
have scaling factors belonging to some special sets (see
Section~\ref{sec:IP-sets}). For an informal pictorial explanation of
our results refer to Section~\ref{sec:informal}.

There are several classes of tilings that have been studied in the
literature (see Section~\ref{sec:definitions}). Roughly, a tiling is
an arrangement of tiles (all of which are copies of some prototiles)
that covers the Euclidean space without overlapping. The classes of
tilings differ in whether when generating the tiles out of the
prototiles one allows rotations, and in whether or not the tiling has
a property ``finite local complexity'' that roughly indicates that
the tiles cannot slide.

Even if the results are slightly different in all these cases, the
ideas of the proofs are very similar. In all cases, we define a
tiling space by considering the closure of all the translations of the
given tiling in the appropriate topology (this topology depends on the
considerations alluded to in the previous paragraph).  One then shows
that translations are commuting homeomorphisms of this tiling space
and one can use Birkhoff's multiple recurrence theorem. The statements
of the results are obtained by unraveling the previous definitions of
the tiling spaces and the meaning of convergence in these spaces.

Our proof mirrors Furstenberg's proof of Gallai's theorem using the
Birkhoff multiple recurrence theorem~\cite{Furstenberg}. Gallai's
theorem is a higher dimensional version of Van der Waerden's
theorem~\cite{Rado}.  We recall that Van der Waerden's theorem states
that any partition of the integers into a finite number of sets has
arbitrarily long arithmetic progressions contained in one of the sets
of the partition.

The analogy of tilings with dynamical systems in $\Z^d$ was noted in
\cite{Rudolph89}. The application of ergodic theory to tilings was
surveyed in \cite{RadinMiles} and \cite{RadinMilesBook}.  The use of
recurrence (either topological or measure-theoretic) to obtain
combinatorial or Ramsey theory results is surveyed in
\cite{Furstenberg}.


Let us emphasize that the results of this paper show that all patterns
can be found in any tiling of the Euclidean space.  Sometimes the fact
that arbitrary sets (or graphs) cannot avoid having some amount of
order is called \emph{Ramsey theory}~\cite{RamseyTheory}.


\section{Tiling Spaces and the Tiling Topology}
\label{sec:definitions}

We reproduce here the requisite definitions of tiling spaces and their
topologies. We refer the reader to the excellent survey~
\cite{Robinson} for further exposition and examples. We have augmented
the definitions to allow for pinwheel and more general tilings
following~\cite{Radin}.

\subsection{Tiling Spaces}

A set $D \subset \R^d$ is called a \emph{tile} if it is compact and equal to
the closure of its interior. A \emph{tiling} of $\R^d$ is a collection
of tiles $\{D_i\}$ such that 
\begin{enumerate}
\item $\bigcup_{i} D_i = \R^d$ -- we say that the tiling
  \emph{covers} $\R^d$.

\item for all $i,j$ with $i \neq j$, $D_i^\circ \cap D_j^\circ =
  \varnothing$ -- we say that the tiling \emph{packs} $\R^d$. 
\end{enumerate}
Two tiles are called \emph{equivalent} if one is a translation of the
other. An equivalence class is called a \emph{prototile}. Two tiles
are called \emph{congruent} of they are related by an orientation
preserving isometry. We will call the equivalence class of congruent
tiles a \emph{congruence prototile}.

\begin{defn}
  Let $\mathcal{T}$ be a finite set of distinct prototiles
  (resp. congruence prototiles) in $\R^d$. We call the tiling space
  $X_{\mathcal{T}}$ associated to $\mathcal{T}$ the set of all tilings
  of $\R^d$ by tiles that are elements of the prototiles
  (resp. congruence prototiles) in $\mathcal{T}$.
\end{defn}

In both cases  the translation of a tiling is still a tiling and hence there
is a natural action of $\R^d$ on $X_{\mathcal{T}}$. We shall use
$T_{\vec{v}}$ to denote both the translation by $\vec{v} \in \R^d$ and
its induced action on $X_{\mathcal{T}}$. If $\mathcal{T}$ consists of
congruence prototiles then there is a natural action of the group of
orientation preserving isometries on $X_{\mathcal{T}}$. 

\begin{defn}\label{patches}
  Let $\mathcal{T}$ be a finite set of distinct (congruence)
  prototiles in $\R^d$. A \emph{$\mathcal{T}$-patch} is subset $x'$ of
  a tiling $x$ such that the union of tiles in $x'$ is a connected
  subset of $\R^d$.\footnote{Some authors require that the union be a
    simply connected subset of $\R^d$, but we will not use that.} We
  call the union of tiles in $x'$ the support of the patch $x'$,
  denoted $\supp(x')$. We call two patches equivalent if one is a
  translation of the other and call an equivalence class a
  \emph{protopatch}. We call two patches congruent if they are related
  by an orientation preserving isometry and call an equivalence class
  of congruent patches a \emph{congruence protopatch}. If
  $\mathcal{T}$ consists of prototiles (resp. congruence prototiles)
  then we denote by $\mathcal{T}^{(n)}$ the collection of protopatches
  (resp. congruence protopatches) consisting of $n$-tiles.
\end{defn}

Several elaborations on this definition can be found in the literature
including distinguishing identical prototiles by coloring or
labeling. Our results can be adapted to such finite extensions but we
will not consider them here. 

\subsection{Complexity of tilings} 

We emphasize that according to our definition tilings have either a
finite number of prototiles or a finite number of congruence
prototiles.

In the following subsections, we will consider three situations of
increasing generality.

\subsubsection{Finite Local Complexity under Translation}
\begin{defn}\label{FLC}
  We say that a tiling space has \emph{finite local complexity under
  translation} if $\mathcal{T}$ consists of a finite number of
  prototiles and $\mathcal{T}^{(2)}$ is finite.  That is, the set of
  pairs of adjacent tiles is finite.
\end{defn}

One example of a tiling that does not have finite local complexity is
the tiling of the plane by a single translated square. The tiles are
necessarily located so that there are countably many straight line
that consists of entirely of tile edges, these are called shear
lines. On each side of a shear line the tiling can be translated along
the shear line and these translations can be chosen
independently. Hence for any square along the line there are
uncountably many configuration of nearby tiles.

Often a finite $\mathcal{T}^{(2)}$ is specified in the description of
the tiling space. This describes a set of local matching rules. These
local matching rules are often specified in terms of edge labels or
colors. These local matching rules can be disposed of by suitably
modifying the prototiles, adding little teeth that make them to match
only on very precise locations. Wang tiles, the classical Penrose kite
and dart tiling, and the tiling by Penrose rhombs all have such local
matching rules. In all these cases the prototiles are polygonal. There
are local matching rules that require the polygons to meet full edge
to edge and, in the case of the Penrose tilings, to preclude the
formation of shear lines.

\subsubsection{Finite Local Complexity under the Euclidean Group}
\begin{defn}\label{FLCE}
  We say that a tiling space has \emph{finite local complexity under
    the Euclidean group} if $\mathcal{T}$ consists of a finite
  number of congruence prototiles and $\mathcal{T}^{(2)}$ is finite.
\end{defn}

The pinwheel tessellation of J. H. Conway was shown to be a tiling of
finite local complexity under the full Euclidean group
in~\cite{Radin}.

\subsubsection{Tilings without Finite Local Complexity}

\begin{defn}\label{NoFLC}
  We say that a tiling space does not have \emph{finite local
    complexity} if $\mathcal{T}$ consists of a finite number of
  congruence prototiles and $\mathcal{T}^{(2)}$ is infinite.
\end{defn}

The simplest example of this type of tiling is the tiling by squares,
discussed above.

\subsection{The Tiling Topologies} 
In this section we will discuss several variations of the 
topology of tilings (which are induced by a metric). 
There are several variations in the literature. All of them 
have in common the fact that two tilings are at small distance 
when they are very similar in a large ball about the origin. 

\subsubsection{The General Metric}
\label{sec:general-metric}

We first give a definition of a metric that applies to all three of
our tiling situations. This is the metric used in \cite[Page
70]{RadinMiles} with proofs in \cite{RadinWolff}. 

Given two tilings $x, x' \in X_{\mathcal{T}}$ we define $d(x,x')$ by
\begin{equation}\label{metric}
  d( x,x') = \sup_{n \in \N} \frac{1}{n} d_H \bigl( B_n(\partial x),
  B_n(\partial x') \bigr) 
\end{equation}
where $d_H$ denotes the usual Hausdorff distance between compact
subsets of $\R^d$, $B_n(\partial x) = B_n \cap \partial x$, $B_n:=\{p
\in \R^d: \|p\|\leq n \}$, and $\partial x = \cup_{D \in x} \partial
D$. In other words, two tilings are close if their skeletons are close
on a large ball about the origin.

The following result is proved  in \cite{Radin,RadinWolff}.
\begin{lem}\label{compactness1}
The  metric  makes the tiling space complete and compact. Moreover, the
action of $\R^d$ by translation on the tiling space is continuous.
\end{lem} 

Proving compactness consists of showing that given a sequence, one can
extract a convergent subsequence.  Given one sequence of tilings, if
we consider any ball, we can extract a subsequence so that the tiles
in the ball converge.  Note that the positions of the tiles in the
ball are given by a finite number of real parameters, which lie on a
bounded set. Then, going to a larger ball, we can extract another
sequence that converges in both balls. We can repeat the argument over
an increasing sequence of balls and then perform a diagonal
argument.  It is easy to verify that the resulting object is indeed a
tiling (if it was not, some violation of the conditions would happen
in a finite ball). It is also easy to verify that the diagonal
sequence indeed converges in the metric.


The verification of continuity of translations is straightforward. We
just note than, if two tilings are very similar in a large ball about
the origin, then the two tilings that result from applying a small
translation will also be very similar in a large ball about the
origin.

Note that the proof indicated here does not use local complexity.  It
only uses the fact that the position of all the tiles in a ball is
indicated by a finite number of parameters.

\subsection{Adapted Metrics}
\label{sec:adapted-metric}

Though the previous metric applies equally well to all three
situations we will give separate metrics for each of the three
situations. Each of these metrics is equivalent to the general metric
but gives more geometric information.

\subsubsection{Finite Local Complexity under Translation}
We will define a metric $d_1$ on a tiling space that has finite local
complexity under the action of translation that makes two tilings
close if after a small translation they agree on a large ball about
the origin~\cite{Rudolph88}.

Let $K \subset \R^d$ be compact. We denote by $x[[K]]$ the collection
of patches $x'$ contained in $x$ with the property that
$\supp(x')\supseteq K$. Let $\|\cdot\|$ denote the Euclidean norm on
$\R^d$ and $B_r = \{ p \in \R^d : \|p\| \leq r \}$. 
\begin{defn}
  We define a metric on the tiling space $X_{\mathcal{T}}$ by 
  \begin{multline*}
    d_1(x,y) = \inf \bigl(\{\frac{1}{\sqrt{2}}\} \cup \{0 < r <
    \frac{1}{\sqrt{2}} : \exists x' \in x[[B_{\frac{1}{r}}]],\\
    y' \in y[[B_{\frac{1}{r}}]], \vec{v} \in B_r,\text{ such that }
    T_{\vec{v}}\, x' = y' \} \bigr).
  \end{multline*}
\end{defn}
The proof that this defines a metric may be found in~\cite{Robinson}
and does not depend on finite local complexity. The only property that
is not immediately obvious is the triangle inequality. The bound of
$\frac{1}{\sqrt{2}}$ is used precisely in the verification of the
triangle inequality.  The r\^ole that finite local complexity plays is
summarized by the following theorem of Rudolph~\cite{Rudolph89}.

\begin{lem}\label{completeness}
  Any tiling space $X_{\mathcal{T}}$ endowed with $d_1$ is
  complete. If $X_{\mathcal{T}}$ has finite local complexity under
  translation then $X_{\mathcal{T}}$ endowed with the metric $d_1$ is
  compact. Moreover, the action of $\R^d$ on $X_{\mathcal{T}}$ by
  translation is continuous.
\end{lem}


\subsubsection{Finite Local Complexity under the Euclidean Group}
We will define a metric $d_2$ on a  tiling space with finite local
complexity under the full Euclidean group that makes two tilings
close if after a small isometry they agree on a large ball about
the origin.

We define a metric on the group of direct isometrics of $\R^d$
\begin{equation*}
  \mathcal{E}^d := \{ T \vec{p} = A \vec{p} + \vec{b} : A \in SO(d),
  \vec{b} \in \R^d \} 
\end{equation*}
by 
\begin{equation*}
  d_{\mathcal{E}^d} (A_1 \vec{p} + \vec{b}_1, A_2 \vec{p} + \vec{b}_2) :=
  \max \{ \|A_1 -A_2\| , \| \vec{b}_1 - \vec{b}_2 \| \}. 
\end{equation*}
Using a common abuse of notation we can write 
\begin{equation*}
  B_r( \mathcal{E}^d ) = \{T \in \mathcal{E}^d : d_{\mathcal{E}^d}( T ,
  \Id) < r\} = \{ T \vec{p} = A \vec{p} + \vec{b} : \|A - \Id \|< r,
  \|\vec{b}\| < r \}.
\end{equation*}
Using these notations we define a metric in a fashion similar to
$d_1$. 
\begin{defn}
  We define a metric on the tiling space $X_{\mathcal{T}}$ by 
  \begin{multline*}
    d_2(x,y) = \inf \bigl(\{\frac{1}{\sqrt{2}}\} \cup \{0 < r <
    \frac{1}{\sqrt{2}} : \exists x' \in x[[B_{\frac{1}{r}}]],\\
    y' \in y[[B_{\frac{1}{r}}]], T \in B_r(\mathcal{E}^d),\text{ such that }
    T\, x' = y' \} \bigr).
  \end{multline*}
\end{defn}
The proof that this is a metric follows almost exactly the proof
in~\cite{Robinson} for $d_1$.

Under this definition any tiling space is complete and the action of
translations is continuous. If the tiling space has finite local
complexity under the Euclidean group then the tiling space endowed
with the metric $d_2$ is compact. The argument is the same as that
sketched for $d$ following Lemma~\ref{compactness1}. 

\subsubsection{Tilings without Finite Local Complexity}
If the tiling space does not have finite local complexity then metrics
which focus on motions of the whole tiling will not give compactness
of the tiling space. In this case we define two tilings to be close if
they agree on a large ball about the origin after a small motion of
each individual tile.

\begin{defn}
  We define a metric on the tiling space $X_{\mathcal{T}}$ by 
  \begin{multline*}
    d_3(x,y) = \inf \bigl(\bigl\{\frac{1}{\sqrt{2}}\bigr\} \cup
    \bigl\{0 < r <
    \frac{1}{\sqrt{2}} : \exists \{t_i\}_{i=1}^m \in x[[B_{\frac{1}{r}}]],\\
    \{s_i\}_{i=1}^m \in y[[B_{\frac{1}{r}}]], \{T_i\}_{i=1}^m \subset
    B_r(\mathcal{E}^d),\text{ such that } T_i\, t_i = s_i \text{ for
      $i =1, \dots, m$ } \bigr\}\bigr).
  \end{multline*}
\end{defn}

Under this definition any tiling space is compact and the action of
translations is continuous. The proof that this is a metric follows
almost exactly the proof in~\cite{Robinson} for $d_1$.

\section{Topological Multiple Recurrence and its Application to Tiling
  Theory}

The following is the Multiple Birkhoff Recurrence Theorem for
commuting homeomorphisms~\cite[Proposition 2.5]{Furstenberg}.  The
same result holds for commuting continuous maps and can be obtained
from this one by passing to the natural extension, but for our
applications, we only require the version for homeomorphisms.

\begin{thm}\label{thm:MBR}
  Let $X$ be a compact metric space and $T_1, \dots, T_l$ commuting
  homeomorphisms of $X$. Then there exists a point $x \in X$ and a
  sequence $n_k \rightarrow \infty$ such that $T^{n_k}_i x \rightarrow
  x$ simultaneously for $i = 1, \dots, l$.
\end{thm}

Furstenberg's original application of this theorem was to prove
Gallai's extension of the Van der Waerden's theorem to higher
dimensions~\cite{Rado}.\footnote{Gallai is also known as Gr\"unwald.}

Our proof follows the same scheme as Furstenberg's proof of
Gallai's theorem ~\cite[Theorem 2.7]{Furstenberg} though the
structure of our topological spaces is quite different.  

\begin{thm}[Main Theorem for Tilings with Finite Local Complexity
  under Translation]\label{MainThm1}
  Let $X_{\mathcal{T}}$ be an $\R^d$-tiling space with finite local
  complexity under translation, and let $x \in X_{\mathcal{T}}$ be an
  arbitrary tiling. Given $\epsilon >0$ and a finite subset $F \subset
  \R^d$ there exists an $n \in \N$, and a patch $p$ contained in $x$
  such that
  \begin{enumerate}
  \item $\supp(p)$ contains a ball of radius $\frac{1}{\epsilon}$ (not
    necessarily centered at the origin),  

  \item for each $\vec{u} \in F$ there exists a vector $\vec{c}$ with
    $\|\vec{c}\:\!\|< \epsilon$ such that
    \begin{equation*}
      T_{n \vec{u} + \vec{c}}\, p \subset x. 
    \end{equation*}
  \end{enumerate}
\end{thm}

\begin{proof}
  Consider the set 
  \begin{equation*}
    X = \overline{\{ T_{\vec{v}}\, x : \vec{v} \in \R^d \}}
  \end{equation*}
  where the closure is taken according to the topology induced by the
  metric $d_1$.  $X$ is a translation invariant compact subset of
  $X_{\mathcal{T}}$. Let $F = \{ \vec{u}_1, \dots , \vec{u}_l\}$ and
  consider the $l$ commuting homeomorphisms of $X$ given by $T_i =
  T_{-\vec{u}_i}$.

  By Theorem \ref{thm:MBR} there exists a point $y \in X$ and a
  sequence $n_k \rightarrow \infty$ such that $T_i^{n_k} y \rightarrow
  y$ for $1 \leq i \leq l$. In particular for large enough $n_k$ we
  have $d_1( T_i^{n_k} y, y ) < \epsilon$. Since $y \in X$ is either a
  translation of $x$ or the limit of translations of $x$ we can find
  $\vec{v} \in \R^d$ such that $d_1(T_{-n_k \vec{u}_i- \vec{v}}\, x,
  T_{-\vec{v}}\, x ) < \epsilon$. By definition of the metric there
  exists $x' \in (T_{-n_k \vec{u}_i - \vec{v}}\,
  x)[[B_{\frac{1}{\epsilon}}]]$, $p_i \in (T_{-\vec{v}}\,
  x)[[B_{\frac{1}{\epsilon}}]]$, and a vector $\vec{c}_i \in \R^d$
  with $\|\vec{c}_i\;\!\|< \epsilon$ such that $T_{-\vec{c}_i} x' =
  p_i$. Now consider $p'$ to be the connected component of
  $\cap_{i=1}^l p_i$ that contains the origin. Since each $p_i$ has
  $B_{\frac{1}{\epsilon}} \subset \supp{p_i}$ we see that
  $B_{\frac{1}{\epsilon}} \subset \supp{p'}$. Now we define $p =
  T_{\vec{v}}\, p'$.  By construction $ p \subset T_{-n_k \vec{u}_i -
    \vec{c}_i}\, x$ and thus we get
  \begin{equation*}
    T_{n_k \vec{u_i} + \vec{c_i}}\, p \subset x
  \end{equation*}
  as required. 
\end{proof}

The only things we have used crucially are the compactness of the
space and the continuity of the action by translation. For
completeness we give proofs for the remaining two cases though we
emphasize that the crucial step is the same in all the proofs.

\begin{thm}[Main Theorem for Tilings with Finite Local Complexity
  under the Euclidean Group]\label{MainThm2}
  Let $X_{\mathcal{T}}$ be an $\R^d$-tiling space of finite local
  complexity, and let $x \in X_{\mathcal{T}}$ be an arbitrary
  tiling. Given $\epsilon >0$ and a finite subset $F \subset \R^d$
  there exists an $n \in \N$, a point $\vec{v} \in \R^d$, and a patch
  $p$ contained in $x$ such that
  \begin{enumerate}
  \item $B_{\frac{1}{\epsilon}}(\vec{v}) \subset \supp(p)$,

  \item for each $\vec{u} \in F$ there exists an $S \in
    B_\epsilon(\mathcal{E}^d)$ such that 
    \begin{equation*}
      T_{n \vec{u}+\vec{v}}\, S \, T_{-\vec{v}}\, p \subset x. 
    \end{equation*}
  \end{enumerate}
\end{thm}

\begin{proof}
  Consider the set 
  \begin{equation*}
    X = \overline{\{ T_{\vec{v}}\, x : \vec{v} \in \R^d \}}
  \end{equation*}
  where the closure is taken according to the topology induced by the
  metric $d_2$.  $X$ is a translation invariant compact subset of
  $X_{\mathcal{T}}$. Let $F = \{ \vec{u}_1, \dots , \vec{u}_l\}$ and
  consider the $l$ commuting homeomorphisms of $X$ given by $T_i =
  T_{-\vec{u}_i}$.

  By Theorem \ref{thm:MBR} there exists a point $y \in X$ and a
  sequence $n_k \rightarrow \infty$ such that $T_i^{n_k} y \rightarrow
  y$ for $1 \leq i \leq l$. In particular for large enough $n_k$ we
  have $d_2( T_i^{n_k} y, y ) < \epsilon$. Since $y \in X$ is either a
  translation of $x$ or the limit of translations of $x$ we can find
  $\vec{v} \in \R^d$ such that $d_2(T_{-n_k \vec{u}_i- \vec{v}}\, x,
  T_{-\vec{v}}\, x ) < \epsilon$. By definition of the metric there
  exists $x' \in (T_{-n_k \vec{u}_i - \vec{v}}\,
  x)[[B_{\frac{1}{\epsilon}}]]$, $p_i \in (T_{-\vec{v}}\,
  x)[[B_{\frac{1}{\epsilon}}]]$, and an isometry $S_i \in
  B_\epsilon(\mathcal{E}^d)$  such
  that $S_i^{-1} x' = p_i$. Now consider $p'$ to be the
  connected component of $\cap_{i=1}^l p_i$ that contains the
  origin. Since each $p_i$ has $B_{\frac{1}{\epsilon}} \subset
  \supp{p_i}$ we see that $B_{\frac{1}{\epsilon}} \subset
  \supp{p'}$. Now we define $p = T_{\vec{v}}\, p'$.  By construction 
  \begin{equation*}
    S_i \, T_{-\vec{v}} \, p = S_i \, p' \subset S_i \, p_i = x'  
  \end{equation*}
  and thus we get
  \begin{equation*}
    T_{n_k \vec{u_i} + \vec{v}}\, S_i \, T_{-\vec{v}} \, p  \subset
    T_{n_k \vec{u_i} + \vec{v}} \, x' \subset x
  \end{equation*}
  as required. 
\end{proof}

\begin{thm}[Main Theorem for Tilings without Finite Local
  Complexity]\label{MainThm3} 
  Let $X_{\mathcal{T}}$ be an $\R^d$-tiling space with $\mathcal{T}$
  finite and let $x \in X_{\mathcal{T}}$ be an arbitrary tiling. Given
  $\epsilon >0$ and a finite subset $F \subset \R^d$ there exists an
  $n \in \N$, $\vec{v} \in \R^d$, and a patch $p=\{t_i\}_{i=1}^m$
  contained in $x$ such that
  \begin{enumerate}
  \item $B_{\frac{1}{\epsilon}}(\vec{v}) \subset \supp(p)$,

  \item for each $\vec{u} \in F$ and $1 \leq i \leq m$ there exists $S_i \in
    B_\epsilon(\mathcal{E}^d)$ such that 
    \begin{equation*}
      T_{n \vec{u}+ \vec{v}}\, S_i \, T_{-\vec{v}} \, t_i \in x. 
    \end{equation*}
  \end{enumerate}
\end{thm}

\begin{proof}
  Consider the set 
  \begin{equation*}
    X = \overline{\{ T_{\vec{v}}\, x : \vec{v} \in \R^d \}}
  \end{equation*}
  where the closure is taken according to the topology induced by the
  metric $d_3$.  $X$ is a translation invariant compact subset of
  $X_{\mathcal{T}}$. Let $F = \{ \vec{u}_1, \dots , \vec{u}_l\}$ and
  consider the $l$ commuting homeomorphisms of $X$ given by $T_i =
  T_{-\vec{u}_i}$.

  By Theorem \ref{thm:MBR} there exists a point $y \in X$ and a
  sequence $n_k \rightarrow \infty$ such that $T_i^{n_k} y \rightarrow
  y$ for $1 \leq i \leq l$. In particular for large enough $n_k$ we
  have $d_3( T_i^{n_k} y, y ) < \epsilon$. Since $y \in X$ is either a
  translation of $x$ or the limit of translations of $x$ we can find
  $\vec{v} \in \R^d$ such that $d_3(T_{-n_k \vec{u}_i- \vec{v}}\, x,
  T_{-\vec{v}}\, x ) < \epsilon$. By definition of the metric there
  exists $x'=\{s_{i,j}\}_{j=1}^{m_i} \in (T_{-n_k \vec{u}_i -
    \vec{v}}\, x)[[B_{\frac{1}{\epsilon}}]]$, $p_i=\{t_{j}\}_{j \in
    J_i} \in (T_{-\vec{v}}\, x)[[B_{\frac{1}{\epsilon}}]]$, and
  isometries $\{S_{i,j}\}_{j \in J_i} \in B_\epsilon(\mathcal{E}^d)$
  such that $S_i^{-1} s_{i,j} = t_{j}$. Now consider $p'$ to be the
  connected component of $\cap_{i=1}^l p_i$ that contains the
  origin. Thus $p' = \{t_j\}_{j \in J}$ for some $J \subset
  \cap_{i=1}^l J_i$. Since each $p_i$ has $B_{\frac{1}{\epsilon}}
  \subset \supp{p_i}$ we see that $B_{\frac{1}{\epsilon}} \subset
  \supp{p'}$. Now we define $p = T_{\vec{v}}\, p'= \{\tau_j\}_{j \in
    J} \subset x$.  By construction
  \begin{equation*}
    S_{i,j} \, T_{-\vec{v}} \, \tau_j = S_{i,j} \, t_j = s_{i,j}   
  \end{equation*}
  and thus we get
  \begin{equation*}
    T_{n_k \vec{u_i} + \vec{v}}\, S_i \, T_{-\vec{v}} \, \tau_j  =
    T_{n_k \vec{u_i} + \vec{v}} \, s_{i,j} \in x
  \end{equation*}
  as required. 
\end{proof}

\subsection{IP-Sets and Dilation Factors }
\label{sec:IP-sets}

There is a refinement of our main theorems that allows some control
over the dilation factors that appear. Unfortunately it does not give
any information on the size of the dilation required. 

\begin{defn}
  A set of positive integers $R$ is called an IP-set if there exists a
  sequence $(p_i)_{i=1}^\infty$ of natural numbers such that $R$
  consists of the numbers $p_i$ together with all finite sums 
  \begin{equation*}
    p_{i_1} + p_{i_2} + \dots + p_{i_k}
  \end{equation*}
  with $i_1 < i_2 < \dots < i_k$. 
\end{defn}

IP-sets appear naturally in situations where recurrence plays a
central r\^ole~\cite{FurstenbergWeiss}~\cite{Bergelson}.  The dilation
factor $n$ in our proofs arises from an application of the Birkhoff
Multiple Recurrence Theorem, Theorem~\ref{thm:MBR}. It is shown
in~\cite[Theorem 2.18]{Furstenberg}, that one can restrict the sets of
numbers that appear in the conclusion of the Birkhoff Multiple
Recurrence Theorem~\ref{thm:MBR} to an a priori given IP-set. Hence
Theorem \ref{MainThm1}, Theorem \ref{MainThm2}, and Theorem
\ref{MainThm3} hold true when we restrict the number $n$ to lie in
some a priori specified IP-set without any modification of the proofs.

\section{An informal pictorial illustrations of 
the results} 
\label{sec:informal} 
For example, given a tiling with finite local complexity (of either
type), a finite set of points in $\R^n$, e.g.,
\begin{center}
\begin{picture}(200,100)
\put(1,50){$1$}
\put(20,50){\vector(0,3){40}}
\put(20,50){\vector(0,-3){40}}
\put(40,30){$\bullet$}
\put(40,40){$\bullet$}
\put(40,50){$\bullet$}
\put(40,60){$\bullet$}
\put(50,30){$\bullet$}
\put(60,30){$\bullet$}
\put(60,40){$\bullet$}
\put(60,50){$\bullet$}
\put(60,60){$\bullet$}
\put(90,30){$\bullet$}
\put(90,40){$\bullet$}
\put(90,50){$\bullet$}
\put(80,60){$\bullet$}
\put(90,60){$\bullet$}
\put(100,60){$\bullet$}
\end{picture}
\end{center}
\noindent 
a size R, and a $\varepsilon >0$. We can find a patch containing a
ball of radius $R$
\vspace{\baselineskip}

\centerline{\epsfysize=1.5truein\epsfbox{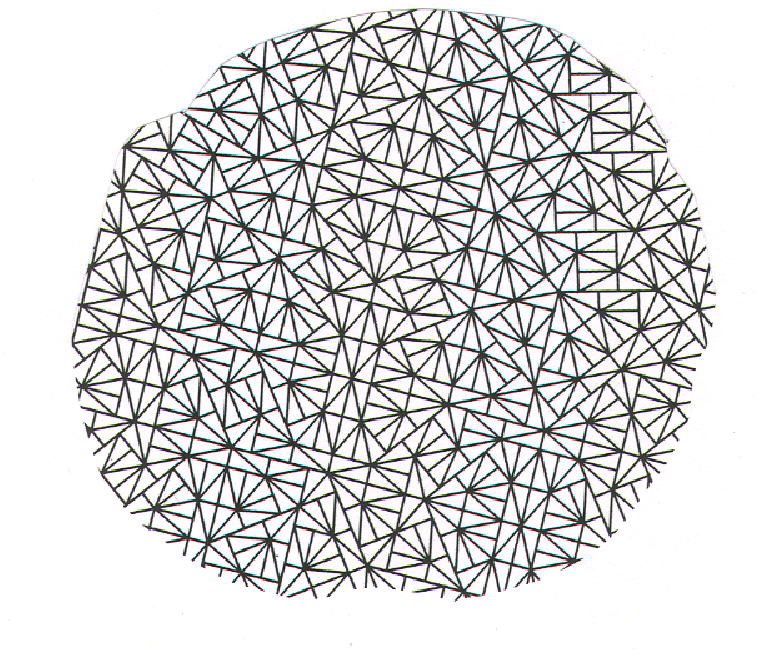}}
and a number $N$, such that the configuration
\begin{center}
\begin{picture}(200,110)
\put(1,50){$N$}
\put(20,50){\vector(0,3){50}}
\put(20,50){\vector(0,-3){40}}
\put(40,20){\epsfysize=.2truein\epsfbox{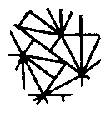}}
\put(40,40){\epsfysize=.2truein\epsfbox{patch-small.eps}}
\put(40,60){\epsfysize=.2truein\epsfbox{patch-small.eps}}
\put(40,80){\epsfysize=.2truein\epsfbox{patch-small.eps}}
\put(60,20){\epsfysize=.2truein\epsfbox{patch-small.eps}}
\put(80,20){\epsfysize=.2truein\epsfbox{patch-small.eps}}
\put(80,40){\epsfysize=.2truein\epsfbox{patch-small.eps}}
\put(80,60){\epsfysize=.2truein\epsfbox{patch-small.eps}}
\put(80,80){\epsfysize=.2truein\epsfbox{patch-small.eps}}
\put(140,20){\epsfysize=.2truein\epsfbox{patch-small.eps}}
\put(140,40){\epsfysize=.2truein\epsfbox{patch-small.eps}}
\put(140,60){\epsfysize=.2truein\epsfbox{patch-small.eps}}
\put(120,80){\epsfysize=.2truein\epsfbox{patch-small.eps}}
\put(140,80){\epsfysize=.2truein\epsfbox{patch-small.eps}}
\put(160,80){\epsfysize=.2truein\epsfbox{patch-small.eps}}
\end{picture}
\end{center}
\noindent 
appears somewhere in the tiling up to an isometry of size less than
$\varepsilon$. The appearance of this configuration may be very far
from the origin.

\section*{Acknowledgments} 
We thank C. Radin, E. A. Robinson and L. Sadun for discussions and
encouragement. M. Combs prepared the figures.  The work of R.L. has
been supported by NSF grants. We thank the referees whose helpful
suggestions resulted in a much better paper. 

\bibliographystyle{alpha}
\bibliography{RecurrenceTiling}

\end{document}